\documentclass[amssymb,amstex,amsmath,10pt,a4paper]{amsart}
\usepackage{amsmath,amssymb,latexsym,amsfonts,amsthm,longtable,url}
\usepackage{epsfig}
\newtheorem{thm}{Theorem}[section]

\newtheorem{lem}[thm]{Lemma}

\newtheorem{rem}[thm]{Remark}

\newcommand{\dist}{{\rm dist}}

\newcommand{\rea}{{\rm Re}}

\newcommand{\Div}{{\rm div}}

\numberwithin{equation}{section}

\begin{document}
\title[New complete embedded minimal surfaces in $\Bbb H\sp2\times\Bbb R$]{New complete embedded minimal surfaces in $\Bbb H\sp2\times\Bbb R$ }
\author[J. Pyo]{Juncheol Pyo}

\begin{abstract}
We construct three kinds of complete embedded minimal surfaces in
$\Bbb H^2\times \Bbb R$. The first is a simply connected, singly
periodic, infinite total curvature surface. The second is an annular
finite total curvature surface. These two are conjugate surfaces
just as the helicoid and the catenoid are in $\mathbb R^3$. The
third one is a finite total curvature
surface which is conformal to $\mathbb S^2\setminus\{p_1,...,p_k\}, k\geq3.$\\

\end{abstract}

\maketitle

%%%%%%%%%%%%%%%%%%%%%%%%%%%%%%%%%%%%%%%%%%%%%%%%%%%%%%%%%%%%%%%%%%%%%%
\textbf{\section{Introduction}}
%%%%%%%%%%%%%%%%%%%%%%%%%%%%%%%%%%%%%%%%%%%%%%%%%%%%%%%%%%%%%%%%%%%%%%
\noindent During recent years the theory of minimal surfaces in $\Bbb
H^2\times \Bbb R$ has been rapidly developed by many mathematicians.
They found some interesting complete minimal surfaces as follows:
the catenoid that is a surface of revolution about the $\Bbb
R$-axis; the helicoid that is ruled by the horizontal geodesic; the
Riemann type minimal surface that is foliated by horizontal circles
and lines; the Scherk type minimal surface that is a minimal graph
over an ideal polygon and is asymptotic to vertical planes (see
\cite{h},\cite{nr},\cite{sa},\cite{st1}).

By Hauswirth and Rosenberg \cite{hr} some properties of complete minimal surfaces of finite
total curvature in $\Bbb H^2\times \Bbb R$ have been revealed. The
vertical plane $\Gamma\times\Bbb R$, where $\Gamma$ is a complete
geodesic in $\Bbb H^2$, is clearly a complete minimal surface of
finite total curvature. Apart from the vertical plane, the only such
surface known to exist is the Scherk type minimal surface. Both
surfaces are simply connected. So Hauswirth and Rosenberg \cite{hr}
raised a natural question: is there a nonsimply connected complete
minimal surface of finite total curvature in $\Bbb H^2\times \Bbb
R$? In particular, is there a minimal annulus of total curvature
$-4\pi$? Note that the rotational catenoid has infinite total
curvature.

In this paper, using the conjugate surface method in $\Bbb H^2\times
\Bbb R$ we construct complete embedded minimal surfaces with $k$
vertical planar ends and total curvature $-4(k-1)\pi$, giving an
affirmative answer to Hauswirth and Rosenberg's question.

The conjugate surface construction is initiated by Smyth \cite{s},
who constructed an embedded minimal disk in a tetrahedron
$T\subset\Bbb R^3$ which is perpendicular to $\partial T$. Then
Karcher, Rossman and others made some complete minimal surfaces by
adopting this conjugate surface method (see
\cite{k},\cite{kps},\cite{r}).

Recently Daniel \cite{d} and Hauswirth, Sa Earp and Toubiana
\cite{hst} generalized the notion of conjugate surface to $\Bbb
H^2\times \Bbb R$. Our construction of new minimal surfaces in $\Bbb
H^2\times \Bbb R$ is based on their theory. We construct a minimal
graph $\Delta^k$ over an infinite triangle in $\Bbb H^2$ such that
$\Delta^k$ is asymptotic to a vertical plane and is bounded by two
horizontal geodesics (one finite, the other infinite) making an
angle of $\pi/k$ and one vertical infinite geodesic (see Figure 1).
It turns out that the conjugate surface of $\Delta^k$ is also a
minimal graph which is perpendicular along its boundary to a horizontal plane and the
two vertical planes making an angle of $\pi/k$ in $\Bbb H^2\times
\Bbb R$.

By reflecting the conjugate surface across these planes we can
construct  a nonsimply connected, genus zero, complete embedded
minimal surface $\Sigma_k$ with total curvature $-4(k-1)\pi$ which
is asymptotic to $k$ vertical planes, $k>1$ (Theorem \ref{eq:t2}).
This is similar to the $k$-noid of $\Bbb R^3$, but a remarkable
difference is that $\Sigma_k$ is embedded in $\Bbb H^2\times \Bbb R$
whereas the $k$-noid has self intersection in $\Bbb R^3$ if
$k\geq3$.

If we extend the minimal graph $\Delta^2$ by $180^\circ$-rotations
about the horizontal boundary geodesics, we obtain a minimal graph
$\Delta_v$ which is bounded by two vertical geodesics (see Figure
2). Rotating $\Delta_v$ by $180^\circ$ about the vertical boundary
geodesics repeatedly, we obtain a simply connected complete embedded
minimal surface which is singly periodic. This surface is different
from  the ruled helicoid of $\Bbb H^2\times \Bbb R$ because it is
not ruled and  because its fundamental piece has finite
total curvature $-4\pi$ whereas the fundamental piece of the ruled
helicoid has infinite total curvature (see Theorem \ref{eq:t3}).
%%%%%%%%%%%%%%%%%%%%%%%%%%%%%%%%%%%%%%%%%%%%%%%%%%%%%%%%%%%%%%%%%%%%%%
\textbf{\section{Preliminaries}}
%%%%%%%%%%%%%%%%%%%%%%%%%%%%%%%%%%%%%%%%%%%%%%%%%%%%%%%%%%%%%%%%%%%%%%
\noindent In $\Bbb H\sp2\times\Bbb R$ we consider the disk model for the
hyperbolic plane $\Bbb H^{2}$ and solid cylinder model for whole
space. Let $x, y$ denote the coordinates in $\Bbb H^{2}$ and $t$
denote the coordinate in $\Bbb R$. Let $\Omega\subset \Bbb
H^{2}\times\{0\}$ be a domain. In $\overline{\Bbb
H}^{2}\times\{0\},$ we denote
$\partial\overline{\Omega}=\partial\Omega\cup\partial_{\infty}\Omega,$
where the boundary part is $\partial\Omega\subset \Bbb
H^{2}\times\{0\}$ and the ideal boundary part is
$\partial_{\infty}\Omega\subset\partial_{\infty}\Bbb
H^{2}\times\{0\}.$ Consider a $\mathcal{C}^{2}$ function $t=u(x,y).$ The
vertical minimal surface equation in $\Bbb H\sp2\times\Bbb R$ is the
following:
\begin{equation}\label{eq:n1}
\Div_{\Bbb H}\Big(\frac{\nabla_{\Bbb H}u}{W_{u}}\Big)=0,
\end{equation}
where $\Div_{\Bbb H}$ and $\nabla_{\Bbb H}$ are the hyperbolic
divergence and gradient respectively and
$W_{u}=\sqrt{1+|\nabla_{\Bbb H}u|^2_{\Bbb H}},$ $|\cdot|_{\Bbb H}$
being the norm in $\Bbb H^{2}.$

In the disk model for $\Bbb H^2$, $$\Bbb H^2=\{(x,y)\in\Bbb R^2|x^2 +y^2<1\},$$ with the metric $ds^2=\left(\frac{2}{1-x^2-y^2}\right)^2(dx^2+dy^2)$, the vertical minimal surface equation (\ref{eq:n1}) becomes as follows:
\begin{eqnarray*}
&&\left(1+D^2(x,y)u_{y}^2\right)u_{xx}+\left(1+D^2(x,y)u_{x}^2\right)u_{yy}-2D^2(x,y)u_{x}u_{y}u_{xy}\\
&&+D(x,y)(xu_{x}+yu_{y})(u_{x}^2+u_{y}^2)=0,
\end{eqnarray*}
where $D(x,y)=\frac{1-x^2-y^2}{2}$.

\indent We refer to the existence theorem of minimal surfaces.\\
\begin{thm}\label{eq:p1}(Corollary 4.1 of \cite{st2})\\
Let $\Omega\subset \Bbb H^{2}\times\{0\}$ be a domain and let $g :
\partial\Omega\cup\partial_{\infty}\Omega\rightarrow \Bbb R$ be a bounded
function everywhere continuous except perhaps at a finite set
$S\subset \partial\Omega\cup\partial_{\infty}\Omega.$ Assume that
the finite boundary $\partial\Omega$ is convex. Then $g$ admits an
extension $u : \overline{\Omega}\setminus S \rightarrow \Bbb R$
satisfying the vertical minimal surface equation (\ref{eq:n1}).
Furthermore, the total boundary of the graph of $u$ (that is the
finite and ideal boundary) is the union of the graph of $g$ on
$(\partial\Omega\cup\partial_{\infty}\Omega \setminus S)$ with the
vertical segments $$\{(q,t)|t \in[A:=\lim_{x\rightarrow
q}\inf_{x\neq q}g(x),~ B:=\lim_{x\rightarrow q}\sup_{x\neq q}g(x)],
~x\in
\partial\Omega\cup\partial_{\infty}\Omega \}$$ at any $q\in S$.
\end{thm}
\begin{thm}\label{eq:c2}(Monotone convergence theorem of \cite{cr})\\
Let $\{u_{n}\}$ be a monotone sequence of solutions of (\ref{eq:n1})
in $\Omega$. If the sequence $\{|u_{n}|\}$ is bounded at one point
of $\Omega$, then there is a non-empty open set $U\subset\Omega$
(the convergence set) such that $\{u_{n}\}$ converges to a solution
of (\ref{eq:n1}) in $U$. The convergence is uniform on compact
subsets of $U$ and the divergence is uniform on compact subsets of
$\Omega-U=V$. $V$ is called the divergence set.
\end{thm}
The following well-known theorems are the maximum principle for
minimal surfaces. It is a special case of a lemma by Schoen
 \cite{sch}, and is proven there.
 \begin{thm}\label{eq:m1}(Maximum principle)
 \begin{itemize}
\item[(1)] (Interior maximum principle) Let
 $\Sigma_{1}$ and $\Sigma_{2}$ be minimal surfaces in $\Bbb H\sp2\times\Bbb
 R$. Suppose $p$ is an interior point of both $\Sigma_{1}$ and
 $\Sigma_{2}$, and suppose $T_{p}(\Sigma_{1})=T_{p}(\Sigma_{2})$. If
 $\Sigma_{1}$ lies on one side of $\Sigma_{2}$ near $p$, then
 $\Sigma_{1}=\Sigma_{2}$;
\item[(2)] (Boundary point maximum principle) Suppose $\Sigma_{1},
 \Sigma_{2}$ have $\mathcal{C}^{2}$-boundaries $C_{1}, C_{2}$. Furthermore,
 suppose the tangent planes of both $\Sigma_{1},
 \Sigma_{2}$ and $C_{1}, C_{2}$ agree at $p$, i.e. suppose
 $T_{p}(\Sigma_{1})=T_{p}(\Sigma_{2})$, $T_{p}(C_{1})=T_{p}(C_{2}).$
 If, near $p$, $\Sigma_{1}$ lies to one side of $\Sigma_{2}$, then
 $\Sigma_{1}=\Sigma_{2}.$
\end{itemize}
\end{thm}
By Daniel \cite{d} and by Hauswirth, Sa Earp and Toubiana
\cite{hst}, we have the following two equivalent concept of
associate and conjugate surfaces.  Let $\Sigma\subset \Bbb
H\sp2\times\Bbb R$ be a surface
 equipped with a connection
$\nabla$. Let $N$ denote its unit normal vector field, $J$ denote
the rotation by angle $\frac{\pi}{2}$ on $T\Sigma$ and $S$ denote a
field of symmetric operator $S_{y}:T_{y}\Sigma\rightarrow
T_{y}\Sigma$ for each $y\in\Sigma.$ Let $T$ be the projection of the
vertical vector $\frac{\partial}{\partial t}$ onto the tangent space
$T\Sigma$ of $\Sigma$ and  $\nu=\langle N,\frac{\partial}{\partial
t}\rangle.$ We have $|T|^{2}+{\nu}^2=1$. Let $TC(\Sigma)$ denote the
total curvature of $\Sigma,$ $TC(\Sigma)=\int_{\Sigma}KdA$ where
$K(p)=\det{S_{p}-(1-|T_{p}|^2)}$ (see, for instance, Daniel \cite{d}
or Hauswirth and Rosenberg \cite{hr} for details). We set
$$S_{\theta}=e^{\theta J}S=(\cos\theta)S+(\sin\theta)JS,$$
$$T_{\theta}=e^{\theta J}T=(\cos\theta)T+(\sin\theta)JT.$$
\begin{thm}\label{eq:p3}(Conjugate minimal surface I, \cite{d})\\
 Let $\Sigma$ be a simply connected surface and
$X:\Sigma\rightarrow \Bbb H^{2}\times \Bbb R$ a conformal minimal
immersion. Let $N$ be the normal, $S$ be the symmetric operator on
$\Sigma$ induced by the shape operator of $X(\Sigma)$.
Let $T$ and $\nu$ be defined as above.\\
Let $z_{0}\in \Sigma.$ Then there exists a unique family
$(X_{\theta})_{\theta\in \Bbb R}$ of conformal minimal immersions
$X_{\theta}:\Sigma \rightarrow \Bbb H^{2}\times \Bbb R$ such that
\begin{itemize}
\item[(1)] $X_{\theta}(z_0)=X(z_0)$ and $dX_{\theta}(z_0)=dX(z_0),$
\item[(2)] the metrics induced on $\Sigma$ by $X$ and $X_{\theta}$
are the same,
\item[(3)] the symmetric operator on $\Sigma$ induced by the shape
operator of $X_{\theta}$ is $S_{\theta},$
\item[(4)] $\frac{\partial}{\partial t}=dX_{\theta}(T_{\theta})+\nu
N_{\theta},$ where $N_{\theta}$ is the unit normal to $X_{\theta}.$
\end{itemize}
Moreover the family $X_{\theta}$ is continuous with respect to
$\theta$, and $X_{0}=X.$ The family of immersions
$(X_{\theta})_{\theta\in \Bbb R}$ is called the associate family of
the immersion $X$. In particular the immersion $X_{\frac{\pi}{2}}$
is called the conjugate immersion of the immersion $X$.
\end{thm}
Let $X=(\varphi, h):\Sigma\rightarrow \Bbb H^{2}\times \Bbb R$ be a
conformal minimal immersion. Then $\varphi$ is a harmonic map to
$\Bbb H^{2}$ and $h$ is a harmonic function. The Hopf differential
of $\varphi$ is the following holomorphic $2$-form:
$$Q\varphi=4\Big<\frac{\partial\varphi}{\partial z},\frac{\partial\varphi}{\partial z}\Big>dz^{2}.$$
Because of conformality of $X$, $Q\varphi=-4\left(\frac{\partial h}{\partial z}\right)^2 dz^2,$ where $z=x+iy$ is a local coordinate on $\Sigma$ and $h=\pm\rea\int2i\sqrt{Q\varphi}dz$.
\begin{thm}\label{eq:p5}(Conjugate minimal surface II, \cite{d}, \cite{hst})\\
Let $X=(\varphi, h):\Sigma\rightarrow \Bbb H^{2}\times \Bbb R$ be a
conformal minimal immersion, and
$X_{\theta}=(\varphi_{\theta},h_{\theta})$ its associate family of
conformal minimal immersions. In particular the immersion
$X_{\frac{\pi}{2}}$ is called the conjugate immersion of the
immersion $X$. Let $h_{\frac{\pi}{2}}$ be the harmonic conjugate of
$h$. Then we have
$$Q\varphi_{\theta}=e^{-2\sqrt{-1}\theta}Q\varphi,~~h_{\theta}=(\cos\theta)h+(\sin\theta)h_{\frac{\pi}{2}}.$$
\end{thm}
Now, we refer to  Krust's type theorem  for minimal vertical graphs and
associate family of surfaces in $\Bbb H^{2}\times \Bbb R$. We call
that $G$ is a vertical graph in $\Bbb H^{2}\times \Bbb R$ if $G$ is
graph of $g$, where $g: \Omega\subset\Bbb H^{2}\rightarrow \Bbb R$.
\begin{thm}\label{eq:p6}(Krust's type theorem, Theorem 14 of \cite{hst})\\
Let $X(\Omega)$ be a minimal vertical graph on a convex domain
$\Omega\subset \Bbb H^{2}$. Then the associate surface
$X_{\theta}(\Omega)$, $\theta\in \Bbb R$ is also a vertical graph.
\end{thm}
We can extend a minimal surface across its special boundary.
\begin{thm}\label{eq:p7}(Schwarz reflection principle, \cite{ro})\\
Suppose a minimal surface $\Sigma\subset \Bbb H^{2}\times \Bbb R$
containing a curve $\Upsilon$ as its boundary.
\begin{itemize}
\item[(1)]  $\Upsilon$ is a horizontal or vertical geodesic line then
$\Sigma$ can be extended smoothly across $\Upsilon$ by
$180^\circ$-rotation about $\Upsilon$.
\item[(2)] $\Upsilon$ lies in a plane. $\Upsilon$ is a geodesic of $\Sigma$ and it is not a
horizontal or vertical geodesic line, and $\Sigma$ meets orthogonal
to the plane along $\Upsilon$ then $\Sigma$ can be
extended smoothly across $\Upsilon$ by reflection through the
plane containing $\Upsilon.$
\end{itemize}
\end{thm}
%%%%%%%%%%%%%%%%%%%%%%%%%%%%%%%%%%%%%%%%%%%%%%%%%%%%%%%%%%%%%%%%%%%%%%%
\textbf{\section{Simply connected complete embedded minimal
surface}}
%%%%%%%%%%%%%%%%%%%%%%%%%%%%%%%%%%%%%%%%%%%%%%%%%%%%%%%%%%%%%%%%%%%%%%%
\begin{lem}\label{eq:t1}
Let $0<\alpha<1$ and integer $k\geq2$ be given. Let $D$ be a domain
in $\Bbb H^{2}$ with $2k$ vertices $p_{2m-1}=\alpha
e^{{\sqrt{-1}}{\frac{(2m-2)\pi}{k}}}$,
$p_{2m}=e^{{\sqrt{-1}}{\frac{(2m-1)\pi}{k}}}$, $m=1,...,k$ and $2k$
sides $A_{m}$ be a geodesic from $p_{2m-1}$ to $p_{2m}$, $B_{m}$ be
a geodesic from $p_{2m}$ to
$p_{2m +1}$, $m=1,...,k$ and $p_{1}=p_{2k+1}$.\\
\indent Then there exists a unique (up to a vertical translation)
embedded minimal surface $\Sigma(\alpha,k)$ which has vertical geodesic lines $V_{m}$ through
the $p_{2m-1},$ $m=1,...,k$, as its boundary and the surface is of finite absolute total curvature at most $-\int KdA\leq(2k-2)\pi$. More precisely,
$\Sigma(\alpha,k)$ is the graph of a function $u : D\rightarrow \Bbb
R$ with $u|_{A_{m}}=+\infty$ and $u|_{B_{m}}=-\infty$, $m=1,...,k$.
\end{lem}
\begin{figure}
\begin{center}
\includegraphics[height=6cm, width=12cm]{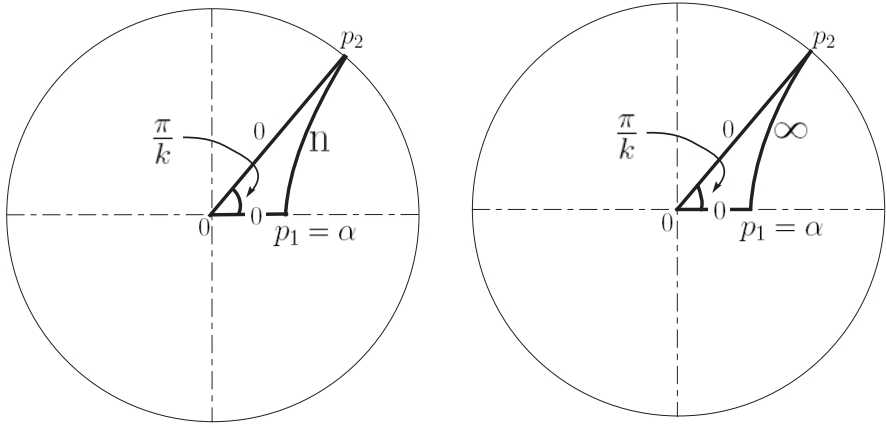}
\end{center}
\caption{Left: the graph of $u_{n}$; Right: the graph of $u$.}
\end{figure}
\begin{proof} Let $L_{1}$ be a geodesic
segment from the origin $0$ of $\Bbb H^{2}$ to $p_{1}$, $L_{2}$ a
geodesic ray from $0$ to $p_{2}$ and $\Gamma$ a geodesic ray from
$p_{1}$ to $p_{2}$. Let $\Omega$ be a convex domain bounded by
$L_{1}$, $L_{2}$ and $\Gamma.$ Let $\widetilde{\Gamma}$ be the complete
geodesic containing $\Gamma$. For each $n\in\mathbb{N}$, let $g$ be a
function on  $\partial\Omega$ such that $g=0$ on $L_{1}\cup L_{2}$
and $g=n$ on $\Gamma.$  By Theorem \ref{eq:p1}, there is a unique
function $u_{n}:\Omega\rightarrow \Bbb R$ satisfying
$u_{n}|_{L_{i}}=0, i=1,2$ and $u_{n}|_{\Gamma}=n$ and the minimal
surface equation (\ref{eq:n1}). By the
maximum principle, $\{u_{n}\}$ is a monotone increasing sequence with respect to $n$.\\
\indent To show that the limit of the sequence $\{u_{n}\}$ exists, we need to
find a suitable barrier.
Let $E$ be the component of $\Bbb H^2\setminus \widetilde{\Gamma}$ which contains the domain $\Omega$.
There exists a function $v\geq0$ defined on $E$, asymptotic to $+\infty$ on $\widetilde{\Gamma}$ and to zero on
$\partial_{\infty}(E)$ and $v$ satisfies the minimal surface equation (\ref{eq:n1}) (see \cite{cr}, \cite{sa}).
 So $v$ is
a suitable barrier for the sequence $\{u_{n}\}$.\\
\indent By the monotone convergence theorem we can find the limit
function $u$ over $\Omega$ of the sequence $\{u_{n}\}$ such that
$u|_{\Gamma}=+\infty$, $u|_{L_{i}}=0$, $i=1,2$ and at $p_{1}$: $\Delta^{k}$, the
graph of $u$, has a vertical geodesic ray as its boundary (see
Figure 1).\\
\indent Since $\Delta^{k}$ lies between the graph of $v$
and the vertical plane $\widetilde{\Gamma}\times\Bbb R$,
outside of a compact part of $\Delta^{k}$ uniformly converges to a vertical plane. We prove in the following that $\Delta^{k}$ has finite total curvature.

More precisely, let $\{q_{j}|q_{j}\in\Gamma\}$ be a sequence such that
$|q_{j}-p_{2}|_{\Bbb R^{2}}$, the Euclidean distance between $q_{j}$
and $p_{2}$ on the disk, is monotone decreasing to zero. Let $\Omega(j)$ be a
compact domain of the convex domain $\Omega$ with $\partial \Omega(j)=\overline{0p}_{1}\cup
\overline{p_{1}q}_{j}\cup\overline{0q}_{j}$, where $\overline{ab}$
indicates the geodesic line segment in $\Bbb H^{2}$ from $a$ to $b$.
For $n,j\in\Bbb N$ fixed, we denote by $u_{n}(j)$ the function on
$\Omega(j)$ which equals $n$ on $\overline{p_{1}q}_{j}$ and zero on
$\overline{0p}_{1}\cup \overline{0q}_{j}$ and satisfies (\ref{eq:n1}). Let $P$ be a point in $\Bbb H^2\times\Bbb R$, we denote $P=(p,t),$ where $p$ is the complex coordinate of $\Bbb H^2$ and $t$ is the coordinate of $\Bbb R$. Denote $v_{1}(j)=(0,0)$, $v_{2}(j)=(p_{1},0)$, $v_{3}(j)=(p_{1},n)$, $v_{4}(j)=(q_{j},n)$ and $v_{5}(j)=(q_{j},0)$, and $\gamma_{i}(j)$ be the geodesic segment in $\Bbb H^2\times\Bbb R$ from $v_{i}(j)$ to $v_{i+1}(j)$, $i=1,\ldots,5$ and $v_{6}(j)=v_{1}(j)$.  The graph of $u_{n}(j)$ denote by $\Delta^{k}_{n}(j)$, is bounded by $\gamma_{i}(j)$, $i=1,\ldots,5$. Applying the Gauss-Bonnet formula to $\Delta^{k}_{n}(j)$, we have:
$$\int_{\Delta^{k}_{n}(j)} K_{n}(j)+\sum^{5}_{i=1}\int_{\gamma_{i}(j)}k_{g,i}(j) +\sum^{5}_{i=1}\theta_{i}(j)=2\pi,$$
where $K_{n}(j)$ is the Gauss curvature function of $\Delta^{k}_{n}(j)$ on the domain $\Omega(j)$ and $0$ on $\Omega-\Omega(j)$, and $k_{g,i}(j)$ is a geodesic curvature function of $\gamma_{i}(j)$, $i=1,\ldots,5$, and $\theta_{i}(j)$ is the exterior angle at $v_{i}(j)$, $i=1,\ldots,5$. By the Gauss equation, the Gauss curvature function $K$ is nonpositive for any minimal surfaces in $\Bbb H\times\Bbb R$ (see \cite{hr}).
Since $k_{g,i}(j)$ is identically zero and $\theta_{1}(j)=\frac{(k-1)}{k}\pi +\angle p_{2}0q_{j}$, $\theta_{i}(j)=\frac{\pi}{2}$, $i=2,\ldots,5$, the total curvature of $\Delta^{k}_{n}(j)$ is $\frac{(1-k)}{k}\pi -\angle p_{2}0q_{j}$. As $j$ goes to infinity, the sequence of $\{u_{n}(j)\}$ converges monotonically to the previous function $u_{n}$ on $\Omega$ and the sequence of $\{\angle p_{2}0q_{j}\}$ converges to zero. By theorem \ref{eq:c2}, the sequence of $\{u_{n}(j)\}$ converges uniformly on compact sets of $\Omega$ to $u_{n} $. By Fatou's lemma, the absolute value of total curvature of $\Delta^{k}_{n}(j)$ is at most $|\frac{(1-k)}{k}\pi|$ for any $n$.

Similarly, as $n$ goes to infinity the absolute value of total curvature of $\Delta^{k}$ is at most $|\frac{(1-k)}{k}\pi|$.

Using the Schwarz reflection principle, we extend
$\Delta^{k}$ about the geodesic
$L_{2}$, extend again about the image of $L_{1}$, again about the
image of $L_{2}$, and so forth. After $2k$ extensions  we get
$\Sigma(\alpha,k)$ which is an
 embedded minimal surface with $2k$ congruent pieces and $k$ vertical geodesics
 passing $p_{2m-1}$, $m=1,...,k$. So the absolute total curvature of $\Sigma(\alpha,k)$ is at most $(2k-2)\pi.$
  By denoting $u$ a extended function defined on $D$ of the previous function $u$,
the proof is completed.
\end{proof}

In case of $k=2$, the $\Delta_{v}=\Sigma(\alpha,2)$ has two vertical
geodesic lines, $V_{1}$ and $V_{2}$. By the Schwarz reflection
principle, we can extend the $\Delta_{v}$ to the complete minimal
surface $\Sigma(\alpha)$ which is singly periodic. Hence the
following theorem holds (see Figure 2).
\begin{figure}
\begin{center}
\includegraphics[height=6cm, width=12cm]{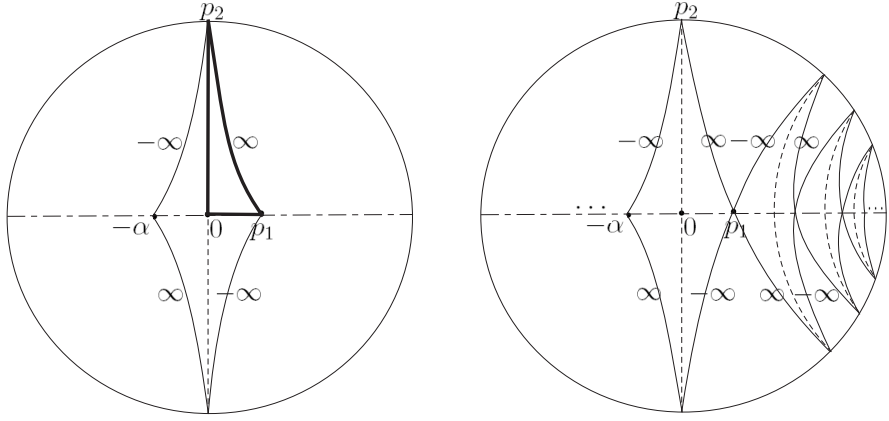}
\end{center}
\caption{Left: the graph of $\Delta_{v}=\Sigma(\alpha,2)$; Right:
Simply-connected complete embedded minimal surface
$\Sigma(\alpha)$.}
\end{figure}
\begin{thm}\label{eq:t3}
In $\Bbb H\sp2\times\Bbb R$, there is a simply connected complete
embedded minimal surface $\Sigma(\alpha)$ which is singly periodic
under the horizontal hyperbolic translation $T_{\alpha}$, $|T_{\alpha}|_{\Bbb
H}=4|\overline{0\alpha}|_{\Bbb H}$. And
the fundamental piece $\Sigma(\alpha)/T_{\alpha}$ has finite total
curvature $-4\pi$.
\end{thm}
\begin{rem}\label{eq:r1}
\begin{itemize}
\item[(1)] We will compute the total curvature of $\Sigma(\alpha)/T_{\alpha}$ in the next section.
\item[(2)] Let $\Delta^{k}=\Delta(\alpha,k)$ be a fundamental piece of
$\Sigma(\alpha,k)$. If $\alpha_{1}\neq\alpha_{2}$,
$\Delta(\alpha_{1},k)$ and $\Delta(\alpha_{2},k)$ can not be
conformally equivalent. So we have one-parameter family of
$\Sigma(\alpha)$ for $\alpha$.
\end{itemize}
\end{rem}
%%%%%%%%%%%%%%%%%%%%%%%%%%%%%%%%%%%%%%%%%%%%%%%%%%%%%%%%%%%%%%%%%%%%%%%
\textbf{\section{Nonsimply connected complete embedded minimal
surface }}
%%%%%%%%%%%%%%%%%%%%%%%%%%%%%%%%%%%%%%%%%%%%%%%%%%%%%%%%%%%%%%%%%%%%%%%
\begin{thm}\label{eq:t2}
For each integer $k\geq2$, there exists  a nonsimply connected
complete embedded minimal surface $\Sigma(k)\subset\Bbb
H\sp2\times\Bbb R$ satisfying the following:
\begin{itemize}
\item[(1)] $\Sigma(k)$ has finite total curvature $-4(k-1)\pi$;
\item[(2)] $\Sigma(k)$ is conformal to a $k$-punctured $2$-dimensional sphere;
\item[(3)] $\Sigma(k)$ is symmetric about $k$ vertical planes and one horizontal
plane.
\end{itemize}
\end{thm}
\begin{proof}
 We take the fundamental piece
$\Delta=\Delta^{k}$ of $\Sigma(\alpha,k)$ in Lemma \ref{eq:t1}, if $k\geq3$ we assume that $\alpha\geq\alpha(k)$, where $\alpha(k)$ is the value that  $\overline{0p}_{1}$ is perpendicular to $\overline{p_{1}p}_{2}$ at $p_{1}$.  The
$\Delta$ is the graph of $u$ over $\Omega$ and is bounded by two
geodesic rays and one geodesic segment. Let $R_{1}$ be a vertical
geodesic ray from $(p_{1},0)$, $R_{2}$ a horizontal geodesic segment
from $(p_{1},0)$ to the origin $(0,0)$ of $\Bbb H^{2}\times\Bbb R$, and $R_{3}$ a
horizontal geodesic ray from $(0,0)$ to $(p_{2},0)$. Here we use the same
notations as in Lemma \ref{eq:t1} except $L_{i}$, $i=1,2$. Since we consider geodesics in $\Bbb H^2\times\Bbb R$, we write $R_{i+1}$, $i=1,2$ instead of $L_{i}$, $i=1,2$.

Let $N$ be a normal vector field
of $\Delta$. Since $\Delta$ is a simply connected, we write its
immersion by $X=(\varphi, h): D^{S}=D\setminus \{S, segment\}\rightarrow
\Delta,$ where $D$ is a closed unit disk in the complex plane $\Bbb
C$. Let $\theta_{1}, \theta_{2} \in (0,2\pi)$ be
 such that $X(e^{\sqrt{-1}\theta_{1}})=(p_{1},0)$ and
$X(e^{\sqrt{-1}\theta_{2}})=(0,0)$.  And let
$c_{1}=\{e^{\sqrt{-1}\theta}|0<\theta\leq\theta_{1}\}$,
$c_{2}=\{e^{\sqrt{-1}\theta}|\theta_{1}\leq\theta\leq\theta_{2}\}$
and $c_{3}=\{e^{\sqrt{-1}\theta}|\theta_{2}\leq\theta<\theta_{3}\}$ such
that
$X(c_{i})=R_{i},$ $i=1,2,3$ and the segment $S$ is $\{e^{\sqrt{-1}\theta}|\theta_{3}\leq\theta\leq2\pi\}$.\\
\indent Let $\Delta_{\frac{\pi}{2}}$ be a conjugate surface of
$\Delta$ with its immersion denoted by
$X_{\frac{\pi}{2}}=(\varphi_{\frac{\pi}{2}}, h_{\frac{\pi}{2}}):
D^{S}\rightarrow \Delta_{\frac{\pi}{2}}$. Let $\gamma_{1}$ be an
arclength parametrization of $R_{1}$. Note that the tangent vector
field $\gamma'_{1}$ is identically $e_{3}$ along $R_{1}$ and
$h_{\frac{\pi}{2}}$ is a harmonic conjugate of $h$. By Theorem
\ref{eq:p5}, on $c_{1}$ we have
$$1=dh\left(\frac{\partial}{\partial s}\right)=dh_{\frac{\pi}{2}}\left(J\frac{\partial}{\partial s}\right).$$
This means that the conormal vector field of
$\Delta_{\frac{\pi}{2}}$ along $\widetilde{R}_{1}$, the conjugate
image of $R_{1}$, is identically $e_{3}.$ So $\widetilde{R}_{1}$
lies on a horizontal plane $\Bbb H^{2}\times \{t\}$ and
$\Delta_{\frac{\pi}{2}}$ is orthogonal to the horizontal
plane along $\widetilde{R}_{1}$. Without loss of generality we assume $t=0$.\\
\indent Let $\gamma_{2}$ be a parametrization of the horizontal
geodesic segment $R_{2}$ and $\widetilde{\gamma}_{2}$ a
parametrization of curve $\widetilde{R}_{2}$, the conjugate image of
$R_{2}$.  Since the conjugate transformation preserves the metric,
$\widetilde{R}_{2}$ is also a geodesic curve on $\Delta_{\frac{\pi}{2}}$. By $(3)$ of Theorem
\ref{eq:p3}, we have
$$0=\langle\overline{\nabla}_{\gamma'_{2}}\gamma'_{2},N\rangle
=\langle S\gamma'_{2},\gamma'_{2}\rangle=\langle
S_{\frac{\pi}{2}}\widetilde{\gamma}'_{2},J\widetilde{\gamma}'_{2}\rangle,$$
where $\overline{\nabla}$ is the connection in $\Bbb H^{2}\times
\Bbb R$. This means that $\widetilde{\gamma}_{2}$ is a line of
curvature.
 On $c_{2}$ we have $$0=dh\left(\frac{\partial}{\partial
s}\right)=dh_{\frac{\pi}{2}}\left(J\frac{\partial}{\partial
s}\right).$$ So the conormal vector field of
$\Delta_{\frac{\pi}{2}}$ along $\widetilde{R}_{2}$ is orthogonal to
$\frac{\partial}{\partial t}$ and $T$, the tangential part of
$\frac{\partial}{\partial t}$. As a result, the curve
$\widetilde{\gamma}_{2}$ is a line of curvature associated to the
field $T$.
\begin{lem}\label{eq:l1}(See the proof of  Proposition 15 of \cite{sot})\\
Let $\Sigma\subset \Bbb H^{2}\times \Bbb R$ be a surface transversal
to each slice $\Bbb H^{2}\times \{t\}$. Let $N$ be a normal field of
$\Sigma$, $T$ a vector field on $\Sigma$ such that $dX(T)$ is the
projection of $\frac{\partial}{\partial t}$ onto tangent plane
$TX(\Sigma)$ and $\nu=\langle N,\frac{\partial}{\partial t}\rangle.$
Let $c:\tau\in I\subset \Bbb R\rightarrow c(\tau)\in \Sigma$ be a
line of curvature associated to the vector field $T.$ Then $c(I)$ is
contained in a vertical totally geodesic plane.
\end{lem}
By Lemma \ref{eq:l1}, $\widetilde{R}_{2}$ is contained in a vertical
plane $\Pi_{1}$. Using an isometry in $\Bbb H^{2}\times\Bbb R$ we
can assume that $\Pi_{1}=\Gamma_{1}\times\Bbb R$, where $\Gamma_{1}$
is a geodesic in $\Bbb H^{2}$ through zero. Let $N_{\frac{\pi}{2}}$
be a normal vector field of $\Delta_{\frac{\pi}{2}}$ and $N_{1}$ the
unit normal vector of $\Pi_{1}$. Since $\widetilde{\gamma}_{2}$ is a
line of curvature and $\widetilde{\gamma}'_{2}\subset T_{\widetilde{\gamma}_{2}}\Pi_{1},$ we
have
$$\frac{d}{ds}\langle N_{\frac{\pi}{2}},N_{1}\rangle
=\langle\overline{\nabla}_{\widetilde{\gamma}'_{2}}N_{\frac{\pi}{2}},N_{1}\rangle
=\langle\mu\widetilde{\gamma}'_{2},N_{1}\rangle=0,$$ where $\mu$ is
a real valued function. So $\langle
N_{\frac{\pi}{2}},N_{1}\rangle=C_{0},$ where $C_{0}$ is a constant.
Since $\nu=\langle N,\frac{\partial}{\partial t}\rangle=1$ at $(0,0)$
and is preserved by the conjugate transformation, $\nu=1$ at
$\widetilde{0}$, the conjugate image of $(0,0)$, $i.e.$
$N_{\frac{\pi}{2}}=\frac{\partial}{\partial t}$.
 So $C_{0}=\langle\frac{\partial}{\partial t},N_{1}\rangle=0.$ Hence
$\Delta_{\frac{\pi}{2}}$ meets $\Pi_{1}$ orthogonally. By the same
argument $\widetilde{R}_{3}$, the conjugate image of $R_{3}$, is
also contained in a vertical plane $\Pi_{2}$ and
$\Delta_{\frac{\pi}{2}}$ meets $\Pi_{2}$ orthogonally. Since
$\widetilde{0}\in \widetilde{R}_{2}$ and $\widetilde{0}\in
\widetilde{R}_{3}$, $\Pi_{1}\cap\Pi_{2}\neq\emptyset$. So we can
assume that $\Pi_{2}=\Gamma_{2}\times\Bbb R$, where
$\Gamma_{2}$ is a geodesic in $\Bbb H^{2}$ through zero. \\
\indent  Because $\nu=1$ at $\widetilde{0}$ and the angle between
$R_{2}$ and $R_{3}$ is $\frac{\pi}{k}$, the angle between
$\widetilde{R}_{2}$ and $\widetilde{R}_{3}$ is also $\frac{\pi}{k}$.
This implies that the angle between $\Pi_{1}$ and $\Pi_{2}$ is
$\frac{\pi}{k}$.  Because $\Delta$ is a graph over the convex domain
$\Omega$, $\Delta_{\frac{\pi}{2}}$ is also a graph over $\Lambda\subset\Bbb
H^{2}\times \{0\}$ by Krust's type theorem. This theorem implies that
$\Delta_{\frac{\pi}{2}}$ is embedded.

\begin{figure}
\begin{center}
\includegraphics[height=8cm, width=6cm]{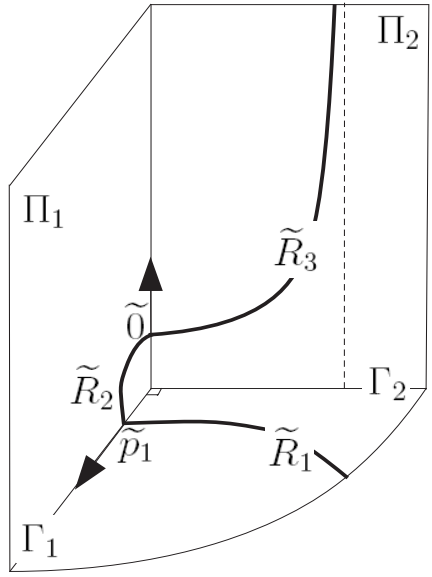}
\end{center}
\caption{In case of $m=2$, the boundary behavior of
$\Delta_{\frac{\pi}{2}}$.}
\end{figure}

We claim that $\Delta_{\frac{\pi}{2}}$ is bounded by $\Pi_{1}, \Pi_{2}$ and $\Bbb H^{2}\times \{0\}$.

We first focus on $\widetilde{R}_{2}$. Define
$d_{{i}}(q)=\dist_{\Bbb R^{2}}(q,\Gamma_{i})$, the Euclidean
distance from $q$ to $\Gamma_{i}$, $i=1,2$  on $\Pi_{1}, \Pi_{2}$.
 Since $\Delta_{\frac{\pi}{2}}$ is a
graph over $\Lambda$, $\widetilde{R}_{2}$ is also graph over
$\Lambda\cap\Gamma_{1}.$ We claim that $d_{1}(q)$ on $\Pi_{1}$
cannot have any interior critical point. Suppose $d_{1}(q)$ has a
local maximum or minimum at $q_{1}\in \widetilde{R}_{2} $.  We have
$\nu(q_{1})=1$. Let $Q_{1}\in R_{2}$ be a preimage of $q_{1}$. Since
$\nu$ is preserved by the conjugate transformation, $\nu(Q_{1})$ is
also one. That is, at $Q_{1}$ the normal vector of $\Delta$ is $e_{3}$. Extend the $\Delta$ along $R_{2}$, the $Q_{1}$ is an interior point of the extended minimal surface. The normal vector of the extended minimal surface at $Q_{1}$ coincides with the one of the horizontal plane $\Bbb H^2\times\{0\}$ and the intersection curve between the extended minimal surface and $\Bbb H^2\times\{0\}$ is just a line. This contradicts to the interior maximum principle.

Since $\nu$ of $\Delta$ varies from $0$ to $1$ as $Q$ varies from $(p_{1},0)$ to $(0,0)$, $\nu$ of $\Delta_{\frac{\pi}{2}}$ varies from $0$ to $1$ as $q$ varies from
$\widetilde{p}_{1}$ to $\widetilde{0}$ along $\widetilde{R}_{2}$. Here $\widetilde{p}_{1}$ is the conjugate image of $(p_{1},0)$.

 Similarly $\widetilde{R}_{3}$ is also a graph over
$\Lambda\cap\Gamma_{2}$ and also $d_{2}(q)$ on $\Pi_{2}$ cannot have
any interior local maximum or minimum. Since the $\nu$ converges to $0$ as
$Q\in R_{3}$ moves toward $(p_{2},0)$, the $\nu$ also converges to
$0$ as $q\in \widetilde{R}_{3}$ varies far away from
$\widetilde{0}$. So the $\widetilde{R}_{3}$ asymptotically
approaches to a vertical geodesic which is orthogonal to $\Gamma_{2}$.

 Now we consider the behavior of $\widetilde{R}_{1}$. The curve $\widetilde{R}_{1}$ cannot intersect with $\Gamma_{1}$ and $\Gamma_{2}$. Suppose not  $\widetilde{R}_{1}$ intersects $\Gamma_{1}$  at $s_{1}$. We extend $\Delta_{\frac{\pi}{2}}$ with respect to $\Pi_{1}$. Then the extended surface $\widetilde{\Delta}_{\frac{\pi}{2}}=\Delta_{\frac{\pi}{2}}\cup{\Delta^{*}}_{\frac{\pi}{2}}$ has a self-intersection, where ${\Delta^{*}}_{\frac{\pi}{2}}$ is the mirror image of $\Delta_{\frac{\pi}{2}}$ with respect to $\Pi_{1}$. But $\widetilde{\Delta}_{\frac{\pi}{2}}$ is the conjugate surface of the graph on the convex domain $\Omega\cup{\Omega^{*}}_{L_{1}}$, where ${\Omega^{*}}_{L_{1}}$ is the $180^{\circ}$-rotated domain about the $L_{1}$. This contradicts to the Krust's type theorem. Similarly, $\widetilde{R}_{1}$ cannot intersect with $\Gamma_{2}$. So $\Delta_{\frac{\pi}{2}}$ is bounded by $\Pi_{1}, \Pi_{2}$ and $\Bbb H^{2}\times \{0\}$.

 We claim that $\widetilde{R}_{1}$ is not convex with respect to $\Lambda$ at any point.
 Suppose $\widetilde{R}_{1}$ is convex at $q_{0}$. Because $\Delta_{\frac{\pi}{2}}$ is a graph over the domain
$\Lambda\subset\Bbb H^{2}\times \{0\}$, near $q_{0}$,
$\Delta_{\frac{\pi}{2}}$ lies on one side of a vertical plane or a
vertical strip of a vertical plane. We extend $\Delta_{\frac{\pi}{2}}$ with respect to $\Bbb H^2\times\{0\}$, then the extended surface intersects with the vertical plane or the
vertical strip of a vertical plane along a point or a line. This contradicts to the
interior maximum principle. Since the $\widetilde{R}_{1}$ is not convex with respect to $\Lambda$ and the length of $\widetilde{R}_{1}$ is infinite, the only option is that $\widetilde{R}_{1}$ goes to the
ideal boundary $\partial_{\infty}\Bbb H^{2}\times\{0\}$ (see Figure 3).

 By the Schwarz reflection principle, we extend
$\Delta_{\frac{\pi}{2}}$ which is of finite absolute total curvature at most $\frac{k-1}{k}\pi$   inductively
 about $\Bbb H^{2}\times \{0\}$, $\Pi_{1}$
and its rotation around $\{0\}\times \Bbb R$-axis by
$\frac{m}{k}\pi$ degrees, $m=1,...,k-1$. In particular, $\Pi_{2}$ is
the rotation of $\Pi_{1}$ around $\{0\}\times \Bbb R$-axis by
$\frac{\pi}{k}$ degrees.  Finally we get $\Sigma(k)$, a complete
embedded surface of finite total curvature with $4k$ congruent fundamental pieces.

First, by Huber's theorem $\Sigma(k)$ is conformal to a  $k$-punctured $2$-dimensional sphere (see \cite{hr} \cite{hu}). This implies that the segment $S$ is nothing but a point. Second, we apply Hauswirth and Rosenberg's curvature estimation \cite{hr} to say that the Hopf map extends meromorphically to each puncture. Moreover, the degree of each pole depends in the number of curves are intersecting horizontal section at infinity. Since this number is one, the degree is zero. So the total curvature of  $\Sigma(k)$ is $-4(k-1)\pi$. Because the fundamental piece of $\Sigma(k)$ is $-\frac{k-1}{k}\pi$, the total curvature of $\Sigma(\alpha)/T_{\alpha}$ is $-4\pi$.
 Third, $\Bbb H^{2}\times
\{0\}$, $\Pi_{1}$ and its rotation around $\{0\}\times \Bbb R$-axis
by $\frac{m}{k}\pi$ degrees, $m=1,...,k-1$ are symmetric planes. By
\cite{hr} each \textit{end}, a conformal parametrization of the
punctured disk, is asymptotic to a vertical plane.
\end{proof}
\begin{rem}
\begin{itemize}
\item[(1)] The $\Sigma(k)$ is similar to the Jorge-Meeks $k$-noid in $\Bbb R^3$
\cite{jm}. And as in Remark \ref{eq:r1}, we have one-parameter
family of minimal surfaces $\Sigma(k)$ with respect to $\alpha$.
\item[(2)] The $\Sigma(\alpha)$ (resp. $\Sigma(2)$) is quite similar to the
Euclidean helicoid (resp. Euclidean catenoid). The $\Sigma(2)$ and
the period $\Sigma(\alpha)/T_{\alpha}$ are conjugate minimal
surfaces, in the sense of Theorem \ref{eq:p3} or Theorem
\ref{eq:p5}.
\end{itemize}
\end{rem}
\textbf{\section{Acknowledgement }} The author would like to
express my gratitude to the referee of this paper, whose remarks
and suggestions have improved the exposition of the text. This
work was supported in part by KRF-2007-313-C00057 and NRF-
2010-0022951.

\vspace{1cm}
\noindent Juncheol Pyo\\
Korea Institute for Advanced Study, 207-43 Cheongryangri 2-dong, Dongdaemun-gu, Seoul 130-722, Korea\\
{\tt e-mail:jcpyo@kias.re.kr}\\

\end{document}